\documentclass[10pt,reqno,twoside,a4paper]{article}
\usepackage {longtable}
\usepackage{amsmath,amssymb,amsfonts}
\usepackage[T2A]{fontenc}
\usepackage[cp866]{inputenc}
\usepackage[english]{babel}
\usepackage{amsthm}
\usepackage[mathscr]{eucal}
\def\{{\protect\lbrace}
\def\}{\protect\rbrace}

\newcommand{\End}{\operatorname{End}}
\newcommand{\Hom}{\operatorname{Hom}}

\newcommand{\Ker}{\operatorname{Ker}}

\newcommand{\Sing}{\operatorname{Sing}}

\begin{document}

\begin{center}
\textbf{Automorphism-Invariant Non-Singular Rings and Modules}
\end{center}
\begin{center}
A.A. Tuganbaev\footnote{National Research University "MPEI", Lomonosov Moscow State University,\\ e-mail: tuganbaev@gmail.com.}
\end{center}

\textbf{Abstract.} A ring $A$ is a right automorphism-invariant right non-singular ring if and only if $A=S\times T$, where $S$ a right injective regular ring and $T$ is a strongly regular ring which contains all invertible elements of its maximal right ring of quotients. Over a ring $A$, each direct sum of automorphism-invariant non-singular right modules is an automorphism-invariant module if and only if the factor ring of the ring $A$ with respect to its right Goldie radical is a semiprime right Goldie ring.

The study is supported by Russian Scientific Foundation (project 16-11-10013).


\textbf{1. Introduction and preliminaries}

All rings are assumed to be associative and with zero identity element; all modules are unitary.
A module $M$ is said to be \textit{automorphism-invariant} if $M$ is invariant under any automorphism of its injective hull. In \cite{DicF69} Dickson and Fuller studied automorphism-invariant modules, when the underlying ring is a finite-dimensional algebra over a field with more than two elements. In \cite[Theorem 16]{ErSS13} Er, Singh and Srivastava proved that a module $M$ is automorphism-invariant if and only if $M$ is a {\it
pseudo-injective} module, i.e.,  for any submodule $X$ of $M$, every monomorphism $X\to M$ can be extended to an endomorphism of the module $M$.
Pseudo-injective modules were studied in several papers; e.g., see \cite{JaiS75}, \cite{Tep75}, \cite{ErSS13}. Automorphism-invariant modules were studied in several papers; e.g., see \cite{AlaEJ05}, \cite{ErSS13}, \cite{GuiS13}, \cite{LeeZ13}, \cite{SinS13}, \cite{Tug13b}, \cite{Tug13d}, \cite{Tug14}, \cite{Tug15}, \cite{Tug17}. 

A ring $A$ is said to be \textit{regular} if every its principal right (left) ideal is generated by an idempotent.
A ring $A$ is said to be \textit{strongly regular} if every its principal right (left) ideal is generated by a central idempotent.
A module $X$ is said to be \textit{injective relative to the module} $Y$ or \textit{$Y$-injective} if for any submodule $Y_1$ of $Y$, every homomorphism
$Y_1\to X$ can be extended to a homomorphism $Y\to X$. A module is said to be \textit{injective} if it is injective with respect to any module.
A module is said to be \textit{square-free} if it does not contain a direct sum of two non-zero isomorphic submodules.
A submodule $Y$ of the module $X$ is said to be \textit{essential} in $X$ if $Y\cap Z\ne 0$ for any non-zero submodule $Z$ of $X$. 
A submodule $Y$ of the module $X$ is said to be \textit{closed} in $X$ if $Y=Y'$ for every submodule $Y'$ of $X$ which is an essential extension of the module $Y$.
We denote by $\Sing X$ the \textit{singular submodule} of the right $A$-module $X$, i.e., $\Sing X$ is a fully invariant submodule of $X$ which consists of all elements $x\in X$ such that $r(x)$ is an essential right ideal of the ring $A$. A module $X$ is said to be \textit{non-singular} if $\Sing X=0$. 

\textbf{Remark 1.1}. In \cite[Theorem 7, Theorem 8, Example 9]{ErSS13} Er, Singh and Srivastava proved the following results.\\
(1)~If $A$ is a right non-singular, right automorphism-invariant ring, then $A=S\times T$, where the ring $S$ is right injective, the module $T_T$ is square-free, any sum of closed right ideals of $T$ is a two-sided ideal which is an automorphism-invariant right $T$-module, and for any prime ideal $P$ of $T$ which is not essential in $T_T$, the factor ring $T/P$ is a division ring.\\
(2)~If $A$ is a right non-singular, right automorphism-invariant prime ring, then the ring $A$ is right injective.\\
(3)~Let $F$ be the field of order 2, $S$ the direct product of a countable set of copies of $F$, and $A = \{(f_n)_{n=1}^{\infty}\in S$: almost all $f_n$ are equal to some $a \in F\}$. Then $A$ is a commutative automorphism-invariant regular ring, but it is not an injective $A$-module.

In connection to Remark 1.1, we will prove Theorem 1.2 which is the first main result of the given paper.

\textbf{Theorem 1.2.} \textit{For a ring $A$, the following conditions are equivalent.}
\begin{enumerate}
\item[\textbf{1)}] 
\textit{$A$ is a right automorphism-invariant right non-singular ring.}
\item[\textbf{2)}] 
\textit{$A$ is a right automorphism-invariant regular ring.}
\item[\textbf{3)}] 
\textit{$A=S\times T$, where $S$ is a right injective regular ring and $T$ is a strongly regular ring which contains all invertible elements of its maximal right ring of quotients.}
\end{enumerate}

\textbf{Remark 1.3}. A module $X$ is said to be {\it quasi-injective} if $X$ is injective relative to $X$, i.e., for any submodule $X_1$ of $X$, every homomorphism $X_1\to X$ can be extended to an endomorphism of the module $X$. Every quasi-injective module is an automorphism-invariant module, since the module $X$ is quasi-injective if and only if $X$ is invariant under any endomorphism of its injective hull; e.g., see \cite[Theorem 6.74]{Lam99}. Every finite cyclic group is a quasi-injective non-injective module over the ring $\mathbb{Z}$ of integers.

For a module $X$, we denote by $G(X)$ the intersection of all submodules $Y$ of the module $X$ such that the factor module $X/Y$ is non-singular. The submodule $G(X)$ is a fully invariant submodule of $X$; it is called the {\it Goldie radical} of the module $X$.

\textbf{Remark 1.4}. In \cite[Theorem 3.8]{KutO80} Kutami and Oshiro proved that any direct sum of non-singular quasi-injective right modules over the ring $A$ is quasi-injective if and only if $A/G(A_A)$ is a semiprime right Goldie ring.

In connection to Remark 1.4, we will prove Theorem 1.5 which is the second main result of the given paper.

\textbf{Theorem 1.5.} \textit{For a ring $A$ with right Goldie radical $G(A_A)$, the following conditions are equivalent.}
\begin{enumerate}
\item[\textbf{1)}] 
\textit{$A/G(A_A)$ is a semiprime right Goldie ring.}
\item[\textbf{2)}] 
\textit{Any direct sum of automorphism-invariant non-singular right $A$-modules is an automorphism-invariant module.}
\item[\textbf{3)}] 
\textit{Any direct sum of automorphism-invariant non-singular right $A$-modules is an injective module.}
\end{enumerate}

The proof of Theorem 1.2 and Theorem 1.5 is decomposed into a series of assertions, some of which are of independent interest. 

We give some necessary definitions. 
A module $X$ is said to be \textit{singular} if $X=\Sing X$. A module $X$ is a \textit{Goldie-radical} module if $X=G(X)$. The relation $G(X)=0$ is equivalent to the property that the module $M$ is non-singular. In the paper, we use well-known properties of $\Sing X$, $G(X)$, non-singular modules and maximal right rings of quotients; e.g., see \cite[Chapter 2]{Goo76}, \cite[Section 7]{Lam99} and \cite[Section 3.3]{Row88}. A module $Q$ is called an \textit{injective hull of the module} $M$ if $Q$ is an injective module and $M$ is an essential submodule of the module $Q$. 
A module $M$ is called a \textit{CS module} if every its closed submodule is a a direct summand of the module $M$. A module $M$ is said to be \textit{uniform} if the intersection of any two non-zero submodules of the module $M$ is not equal to zero. A module $M$ is said to be \textit{finite-dimensional} if $M$ does not contain an infinite direct sum of non-zero submodules.

A ring $A$ is called a \textit{right Goldie ring} if $A$ is a right finite-dimensional ring with the maximum condition on right annihilators. A ring $A$ is said to be \textit{reduced} if $A$ does not have non-zero nilpotent elements. A ring without non-zero nilpotent ideals is said to be \textit{semiprime} ring. A ring $A$ is said to be \textit{right strongly semiprime} \cite{Han75} if any its ideal, which is an essential right ideal, contains a finite subset with zero right annihilator. A ring is said to be \textit{right strongly prime} \cite{Han75b} if every its non-zero ideal contains a finite subset with zero right annihilator. 

\textbf{Remark 1.6.} Every right strongly semiprime ring is a right non-singular semiprime ring \cite{Han75}. It is clear that every right strongly prime ring is right strongly semiprime. The direct product of two finite fields is a finite commutative strongly semiprime ring which is not strongly prime. The direct product of a countable number of fields is an example of a commutative semiprime non-singular ring which is not strongly semiprime. All finite direct products of rings without zero-divisors and all finite direct products of simple rings are right and left strongly semiprime rings. 

\textbf{Remark 1.7.} If $A$ is a semiprime right Goldie ring, then it is well known\footnote{For example, see \cite[Theorem 3.2.14]{Row88}.} that every essential right ideal of the ring $A$ contains a non-zero-divisor. Therefore, all semiprime right Goldie rings are right strongly semiprime. In particular, all right Noetherian semiprime rings are right strongly semiprime.

\textbf{2. Automorphism-Invariant Nonsingular Rings}

\textbf{Lemma 2.1 \cite[Section 3.3.]{Row88}.} \textit{Let $A$ be a right non-singular ring with maximal right ring of quotients $Q$. Then $Q$ is an injective right regular ring and $Q$ can be naturally identified with the ring $\End Q_A$ and $Q_A$ is an injective hull of the module $A_A$.}

\textbf{Lemma 2.2.} \textit{If $A$ is a right non-singular ring with maximal right ring of quotients $Q$, then $A$ is a right automorphism-invariant ring if and only if $A$ contains all invertible elements of the ring $Q$.}

Lemma 2.2 follows from Lemma 2.1.

\textbf{Lemma 2.3 \cite[Chapter 12, 5.1--5.4.]{Stm75}.} \textit{Let $A$ be a reduced ring with maximal right ring of quotients $Q$. Then ring $A$ is right and left non-singular. If each of the closed right ideals of the ring $A$ is an ideal, then $A$ is a reduced ring and $Q$ is a right and left injective strongly regular ring.}

\textbf{Lemma 2.4.} \textit{Let $A$ be a right non-singular ring in which all closed right ideals are ideals. Then the ring $A$ is reduced.}

\textbf{Proof.} Let $a$ be an element of $A$ with $a^2=0$. There exists a closed right ideal $B$ of $A$ such that $B\cap aA=0$ and $B+aA$ is an essential right ideal of $A$. By assumption, the closed right ideal $B$ is an ideal. Therefore, $aB=0$, since $aB$ is contained in the intersection of $aA$ and $B$. Then $a(B+aA)=0$ and $B+aA$ is an essential right ideal. Since $A$ is right non-singular, $a=0$.\hfill$\square$

\textbf{Lemma 2.5.} \textit{If $A$ is a right automorphism-invariant right non-singular ring, then $A=S\times T$, where $S$ is a right injective regular ring and $T$ is a strongly regular ring which contains all invertible elements of its maximal right ring of quotients.}

\textbf{Proof.} By Remark 1.1(1) and Lemma 2.1, $A=S\times T$, where $S$ is a right injective regular ring, $T$ is a right automorphism-invariant right non-singular ring, and any closed right ideal of $T$ is an ideal. By Lemma 2.4, $T$ is a reduced ring. Let $Q$ be the maximal right ring of quotients of the ring $T$. By Lemma 2.3, $T$ is a reduced ring and $Q$ is a right and left injective strongly regular ring. To prove that $T$ is a strongly regular ring, it is sufficient to prove that an arbitrary element $t$ of the ring $T$ is the product of a central idempotent and an invertible element. Since $t$ is an element of the strongly regular ring $Q$, we have that $t=eu$, where $e$ is a central idempotent of the ring $Q$ and $u$ is an invertible element of the ring $Q$. By Lemma 2.2, $T$ contains all invertible elements of the ring $Q$. Therefore, $u\in T$. Then $e=tu^{-1}\in T$ and every element of the ring $T$ is the product of a central idempotent and an invertible element.\hfill$\square$

\textbf{Remark 2.6. The completion of the proof of Theorem 1.2.} In Theorem 1.2, the implication 1)\,$\Rightarrow$\,3) follows from Lemma 2.5, the implication 3)\,$\Rightarrow$\,2) follows from the property that the direct product of regular rings $S$ and $T$ is a regular ring, and the implication 2)\,$\Rightarrow$\,1) follows from the property that every regular ring is right and left non-singular.

\textbf{Corollary 2.7.} \textit{If $A$ is a right automorphism-invariant right non-singular indecomposable ring, then $A$ is a right injective ring}; see Remark 1.1(2).

Corollary 2.7 follows from Theorem 1.2 and the property that every strongly regular indecomposable ring is a division ring and, consequently, a right injective ring.

\textbf{Corollary 2.8.} \textit{Let $A$ be a right automorphism-invariant right non-singular ring which does not contain an infinite set of non-zero central orthogonal idempotents. Then $A$ is a right injective ring.}

Corollary 2.8 follows from Corollary 2.7 and the property that every ring, which does not contain an infinite set of non-zero central orthogonal idempotents, is a finite direct product of indecomposable rings.

\textbf{3. Automorphism-Invariant Non-singular Modules}

{\bf Lemma 3.1.} {\it Let $A$ be a ring and $X$ a right $A$-module which is not an essential extension of a singular module. Then there exists a
non-zero right ideal $B$ of the ring $A$ such that the module $B_A$ is isomorphic to a submodule of the module $X$.}

{\bf Proof.} Since the module $X$ is not an essential extension of a singular module, there exists an element $x\in X$ such that $xA$ is a non-zero non-singular module. Since $xA\cong A_A/r(x)$ and the module $xA$ is non-singular, the right ideal $r(x)$ is not an essential. Therefore, there exists a non-zero right ideal $B$ with $B\cap r(x)=0$. In addition, there exists an epimorphism $f\colon A_A\to xA$ with kernel $r(x)$. Since $B\cap \Ker f=0$, we have that $f$ induces the monomorphism $g\colon B\to xA$. Therefore, $xA$ contains the non-zero submodule $g(B)$ which is isomorphic to the module $B_A$.~\hfill$\square$

\textbf{Lemma 3.2.} \textit{Let $A$ be a ring, $G=G(A_A)$ the right Goldie radical of the ring $A$, $h\colon A\to A/G$ the natural ring epimorphism and $X$ a non-singular non-zero right $A$-module.}
\begin{enumerate}
\item[\textbf{1)}]
\textit{If $B$ is a essential right ideal of the ring $A$, then $h(B)$ is an essential right ideal of the ring $h(A)$.} 
\item[\textbf{2)}] 
\textit{If $B$ is a right ideal of the ring $A$ such that $G\subseteq B$ and $h(B)$ is an essential right ideal of the ring $h(A)$, then $B$ is an essential right ideal of the ring $A$.} 
\item[\textbf{3)}] 
\textit{For any right $A$-module $M$, the module $MG$ is contained in the Goldie radical of $M$.}
\item[\textbf{4)}] 
\textit{$XG=0$ and the natural $h(A)$-module $X$ is non-singular. In addition, if $Y$ is an arbitrary non-singular right $A$-module, then $YG=0$ and the $h(A)$-module homomorphisms $Y\to X$ coincide with the $A$-module homomorphisms $Y\to X$. Therefore, $X$ is an $Y$-injective $A$-module if and only if $X$ is an $Y$-injective $h(A)$-module. The essential submodules of the $h(A)$-module $X$ coincide with the essential submodules of the $A$-module $X$.}
\item[\textbf{5)}] 
\textit{$X$ is an injective $h(A)$-module if and only if $X$ is an injective $A$-module.}
\item[\textbf{6)}] 
\textit{$X_{h(A)}$ is a uniform module $($resp., an essential extension of a direct sum of uniform modules$)$ if and only if $X_A$ is a uniform module $($resp., an essential extension of a direct sum of uniform modules$)$.}
\item[\textbf{7)}] 
\textit{$X_A$ is an essential extension of a direct sum of modules each of them is isomorphic to some non-zero right ideal of the ring $A$.}
\item[\textbf{8)}] 
\textit{If the ring $A$ is right finite-dimensional, then $X_A$ is an essential extension of a direct sum of modules each of them is isomorphic to some non-zero uniform right ideal of the ring $A$.}
\item[\textbf{9)}] 
\textit{If the ring $h(A)$ is right finite-dimensional, then $X_{h(A)}$ is an essential extension of a direct sum of modules each of them is isomorphic to some non-zero uniform right ideal of the ring $h(A)$.}
\end{enumerate}

\textbf{Proof.} \textbf{1.} Let us assume that $h(B)$ is not an essential right ideal of the ring $h(A)$. Then there exists a right ideal $C$ of the ring $A$ such that $C$ properly contains $G$ and $h(B)\cap h(C)=h(0)$. Since $h(B)\cap h(C)=h(0)$, we have that $B\cap C\subseteq G$. Since $C$ properly contains the closed right ideal $G$, we have that $C_A$ contains a non-zero submodule $D$ with $D\cap G=0$. Since $B$ is an essential right ideal, $B\cap D\ne 0$ and $(B\cap D)\cap G=0$. Then $h(0)\ne h(B\cap D)\subseteq h(B)\cap h(C)=h(0)$. This is a contradiction.

\textbf{2.} Let us assume that $B$ is not an essential right ideal of the ring $A$. Then $B\cap C=0$ for some non-zero right ideal $C$ of the ring $A$ and $G\cap C\subseteq B\cap C=0$. Therefore, $h(C)\ne h(0)$. Since $h(B)$ is an essential right ideal of the ring $h(A)$, we have that $h(B)\cap h(C)\ne h(0)$. Let $h(0)\ne h(b)=h(c)\in h(B)\cap h(C)$, where $b\in B$ and $c\in C$. Then $c-b\in G\subseteq B$. Therefore, $c\in B\cap C=0$ and $h(c)=h(0)$. This is a contradiction.

\textbf{3.} For any element $m\in M$, the module $mG_A$ is a Goldie-radical module, since $mG_A$ is a homomorphic image of the Goldie-radical module $G$. Therefore, $mG\subseteq G(M)$ and $MG\subseteq G(M)$.

\textbf{4.} By 3, $XG=0$. Let us assume that $x\in X$ and $xh(B)=0$ for some essential right ideal $h(B)$, where $B=h^{-1}(h(B))$ is the complete pre-image of $h(B)$ in the ring $A$. By 2), $B$ is an essential right ideal of the ring $A$. Then $xB=0$ and $x\in \Sing X=0$. Therefore, $X$ is a non-singular $h(A)$-module. The remaining part of 4 is directly verified.

\textbf{5.} Let $R$ be one of the rings $A$, $h(A)$ and $M$ a right $R$-module. By Lemma 1(4), the module $M$ is injective if and only if $M$ is injective relative to the module $R_R$. Now the assertion follows from 4.

\textbf{6.} The assertion follows from 4.

\textbf{7.} Let $\mathcal{M}$ be the set of all submodules of the module $X$ which are direct sums of modules each of them is isomorphic to a non-zero right ideal of the ring $A$. The set $\mathcal{M}$ is not empty by Lemma 3.1. There exists a partial order in $\mathcal{M}$ such that for any $M,M'\in \mathcal{M}$, the relation $M\lneqq M'$ is equivalent to the property that $M'=M\oplus N$ for some $N\in \mathcal{M}$. By the Zorn lemma, the set $\mathcal{M}$ contains at least one maximal element $K$. 

Let us assume that $K$ is not an essential submodule of the module $X$. Then there exists a non-zero submodule $L$ of the non-singular module $X$ with $K\cap L=0$. By Lemma 3.1, there exists a non-zero right ideal $B$ of the ring $A$ such that the module $B_A$ is isomorphic to some submodule $L'$ of the module $L$. This contradicts to the property that $K$ is a maximal element of the set $\mathcal{M}$.

\textbf{8.} Since the ring $A$ is right finite-dimensional, every non-zero right ideal of the ring $A$ is an essential extension of a finite direct sum of non-zero uniform right ideals. Now the assertion follows from 7.

\textbf{9}. The assertion follows from 6 and 8.~\hfill$\square$

\textbf{Remark 3.3.} Let $M$ be an automorphism-invariant non-singular module. In \cite[Theorem 3, Theorem 6(ii)]{ErSS13} Er, Singh and Srivastava proved that $M = X\oplus Y$, where $X$ is a quasi-injective non-singular module, $Y$ is an automorphism-invariant non-singular square-free module, the modules $X$ and $Y$ are injective relative to each other, $\Hom (X,Y ) = 0 = \Hom (Y,X)$ and $\Hom (D_1,D_2) = 0$ for any two submodules $D_1$, $D_2$ of the module $Y$ with $D_1\cap D_2 = 0$. In addition, for any set $\{K_i\,|\, i\in I\}$ of closed submodules of $Y$, the submodule $\sum _{i\in I}K_i$ is an automorphism-invariant module.

\textbf{Remark 3.4.} Let $M$ be a direct sum of CS modules $M_i$, $i\in I$. In \cite[Corollary 15]{LeeZ13} Lee and Zhou proved that $M$ is a quasi-injective module if and only if $M$ is an automorphism-invariant module.

\textbf{Remark 3.5.} Let $A$ be a ring with right Goldie radical $G(A_A)$.\\ 
In \cite[Theorem 3.4]{KutO80} Kutami and Oshiro proved that the factor ring $A/G(A_A)$ is a right strongly semiprime ring if and only if every non-singular quasi-injective right $A$-module is injective.\\
In \cite[Theorem 3.8]{KutO80} Kutami and Oshiro proved that the factor ring $A/G(A_A)$ is a semiprime right Goldie ring if and only if every direct sum of non-singular quasi-injective right $A$-modules is quasi-injective.

\textbf{Lemma 3.6.} \textit{Let $A$ be a ring with right Goldie radical $G(A_A)$ and $M$ an automorphism-invariant non-singular right $A$-module which is an essential extension of a direct sum of uniform modules.}
\begin{enumerate}
\item[\textbf{1)}]
\textit{$M$ is an essential extension of some quasi-injective non-singular module $K$ which is direct sum of uniform modules closed in $M$.} 
\item[\textbf{2)}] 
\textit{If the factor ring $A/G(A_A)$ is a right strongly semiprime ring, then $M$ is an injective module.}
\end{enumerate}

\textbf{Proof.} \textbf{1.} By Remark 3.3, $M = X\oplus Y$, where $X$ is a quasi-injective module, $Y$ is an automorphism-invariant square-free module. Therefore, we can assume that $M$ is an automorphism-invariant square-free module. Since $M$ is an essential extension direct sum of uniform submodules, $M$ is an essential extension of some module $K$ which is the direct sum of uniform closed submodules $K_i$ of $M$, $i\in I$. By Remark 3.3, $K$ is an automorphism-invariant module. Since every uniform module is a CS module, $K$ is a quasi-injective module by Remark 3.4.

\textbf{2.} By 1, $M$ is an essential extension of some quasi-injective non-singular module $K$. By Remark 3.5, $K$ is an injective essential submodule of the module $M$. Therefore, $K$ is an essential direct summand of the module $M$. Then $M=K$ and $M$ is an injective module.~\hfill$\square$

\textbf{Lemma 3.7.} \textit{Let $A$ be a ring with right Goldie radical $G(A_A)$ and $M$ an automorphism-invariant non-singular right $A$-module. If the factor ring $A/G(A_A)$ is a semiprime right Goldie ring, then $M$ is an injective module.}

\textbf{Proof.} By Lemma 3.2(9), $M$ is an essential extension of a direct sum of uniform modules. In addition, the semiprime right Goldie ring $A/G(A_A)$ is a right strongly semiprime ring \cite{Han75}. By Lemma 3.6(2), $M$ is an injective module.~\hfill$\square$

\textbf{Remark 3.8. The completion of the proof of Theorem 1.5.} In Theorem 1.5, the implications 3)\,$\Rightarrow$\,2)\,$\Rightarrow$\,1) are obvious. 

1)\,$\Rightarrow$\,3). Let $M$ be the direct sum of automorphism-invariant non-singular right $A$-modules $M_i$, $i\in I$. By Lemma 3.7, each of the modules $M_i$ is injective. By Remark 3.5, $M$ is a quasi-injective module. By Remark 3.5, $M$ is an injective module.~\hfill$\square$

\end{document}